\documentclass{amsart}

\usepackage{amsmath, amssymb}
\usepackage{amsfonts}
\usepackage[all]{xy}

\numberwithin{equation}{section}

\DeclareMathOperator{\Aut}{Aut}

\DeclareMathOperator{\id}{id}
\DeclareMathOperator{\Hom}{Hom}
\DeclareMathOperator{\Mod}{mod-}
\DeclareMathOperator{\Br}{Br}
\DeclareMathOperator{\Rep}{Rep}
\DeclareMathOperator{\Arg}{Arg}

\newcommand{\norm}[1]{\lVert #1 \rVert}
\newcommand{\abv}[1]{\lvert #1 \rvert}
\renewcommand{\Im}{\text{Im\,}}
\renewcommand{\Re}{\text{Re\,}}

\def\bs{\backslash}
\def\lto{\longrightarrow}
\def\lam{\lambda}
\def\endproof{\hfill{$\square$} \medskip}
\def\C{\mathbb{C}}
\def\R{\mathbb{R}}

\def\Z{\mathbb{Z}}
\def\alp{\alpha}
\def\sli{\mathcal{P}}
\def\D{\mathcal{D}}
\def\KQ{K(\D_Q)}
\def\DQ{D^{b}(\A_Q)}

\def\H{\mathbb{H}}
\def\Xreg{X_{\mathrm{reg}}}
\def\class{\mathcal{C}}
\def\Stab{\mathrm{Stab}}

\def\op{\mathrm{op}}
\def\rd{\partial}
\def\supp{\mathrm{supp}}
\def\A{\mathcal{A}}
\def\nil{\mathrm{nil}}

\def\CY{\mathrm{CY}}

\newtheorem{thm}{Theorem}[section]
\newtheorem{defi}[thm]{Definition}
\newtheorem{prop}[thm]{Proposition}
\newtheorem{cor}[thm]{Corollary}
\newtheorem{rem}[thm]{Remark}
\newtheorem{lem}[thm]{Lemma}
\newtheorem{conjecture}[thm]{Conjecture}

\begin{document}

\title[Stability conditions for preprojective algebras]
{Stability conditions for preprojective algebras and root systems of Kac-Moody Lie algebras}

\author{Akishi Ikeda}

\begin{abstract}
The aim of this paper is to study the space of  
stability conditions on the bounded derived category of 
nilpotent modules over the preprojective algebra associated with a quiver without loops. 
We describe this space as a covering space of some open set  determined by the root system of the 
Kac-Moody Lie algebra associated with the quiver.  
\end{abstract}
\maketitle

\section{Introduction}
The notion of a stability condition on a triangulated category was introduced 
by T. Bridgeland in \cite{Br1}. He also showed that 
the space of all stability conditions $\Stab(\D)$ 
on a triangulated category $\D$ 
has the structure of a complex manifold 
and there is a local isomorphism
\begin{equation*}
\pi \colon \Stab(\D) \longrightarrow \Hom_{\Z}(K(\D),\C)
\end{equation*}
where $K(\D)$ is the Grothendieck group of $\D$. 

In geometric setting, the spaces of stability conditions for some 
$2$-Calabi-Yau ($\CY_2$) triangulated categories were studied in \cite{Br2,Br3,IUU,Ok2}. 
When $\D$ is the bounded derived category of coherent sheaves of a K$3$ surface \cite{Br2} or 
a certain subtriangulated category of coherent sheaves of a resolution of 
a Kleinian singularity \cite{Br3}, 
it was shown that the distinguished connected component of $\Stab(\D)$ is a covering space 
of some open subset of $\Hom_{\Z}(K(\D),\C)$ determined by 
the data of the root system 
and it has been clarified how the group of deck transformations acts on it. 
Further, the problems of connectedness and simply connectedness 
of $\Stab(\D)$ for 
some particular $\D$ were solved in \cite{IUU,Ok2,ST}. 

Other examples of $\CY_2$ categories come from representation theories. 
It is known that the derived category of 
representations of the preprojective algebra associated 
with a quiver $Q$ is $\CY_2$ category,  
if $Q$ is not of ADE type (see \cite{Kel}). 
$\CY_2$ categories coming from Kleinian singularities 
also have a description in this context, 
as the derived categories of preprojective algebras defined by the corresponding Dynkin diagrams (see \cite{Br3}). 

The aim of this paper is to study the distinguished connected component of 
$\Stab(\D_Q)$ for the bounded derived category $\D_Q$ of the preprojective algebra associated with a quiver $Q$ without loops ($Q$ is not of ADE type). 
We describe this space as a covering space of some open set  determined by the root system of the 
Kac-Moody Lie algebra associated with $Q$.  
This generalizes the result by Thomas for quivers of A type \cite{T}
and Bridgeland for quivers of 
ADE type and affine ADE type \cite{Br3} to arbitrary 
quivers without loops. 
Many basic ideas of this paper come from \cite{Br2,Br3}.

\subsection{Summary of results}
Let $Q$ be a connected finite quiver without loops ($1$-cycles) and let $\{1,\dots,n\}$ 
be the vertices of $Q$. Further, assume that 
the underlying graph $\underline{Q}$, which is a graph obtained 
by forgetting orientations of arrows in $Q$, is not of ADE type. 
For  the quiver $Q$, we can define 
a $\C$-algebra $\Pi(Q)$, called the 
preprojective algebra of $Q$, 
which has a natural grading $\Pi(Q) =\oplus_{i \ge 0}\Pi(Q)_i$ 
by the length of paths.  
A right $\Pi(Q)$-module $M$ is called nilpotent if there is some positive integer $k$ such that 
$M \Pi(Q)_l =0$ for all $l \ge k$. 
Let $\A_Q$ be an abelian category of finite dimensional 
nilpotent right $A$-modules and $\D_Q:=D^b(\A_Q)$ be the 
bounded derived category of $\A_Q$. 
It is known that the triangulated category $\D_Q$ is a $\CY_2$ triangulated category. 

The abelian category $\A_Q$ is finite length with finitely many simple modules 
$S_1,\dots,S_n$ corresponding to vertices of $Q$. 
Hence, every object of $\A_Q$ has the Jordan-H\"older filtration by the simple modules. 
The Grothendieck group of $\D_Q$ is given by 
\begin{equation*}
K(\D_Q) \cong \bigoplus_{i=1}^n \Z [S_i].
\end{equation*}
There is  a natural bilinear form $\chi \colon K(\D_Q) \times K(\D_Q) \to \Z$, called 
the Euler form, defined by 
\begin{equation*}
\chi(E,F):=\sum_{i \in \Z} (-1)^i \dim_{\C} \Hom_{\D_Q}(E,F[i]), 
\quad \text{for } E,F \in \D_Q. 
\end{equation*}
The $\CY_2$ property of $\D_Q$ implies that the Euler form $\chi$ is 
symmetric. 

For $Q$, we associate an $n \times n$ symmetric matrix $A_Q$, called the generalized Cartan matrix (GCM for short), by 
$(A_Q)_{ij} := 2\delta_{ij} - (q_{ij} + q_{ji})$ 
where $q_{ij}$ is the number of arrows from $i$ to $j$ in $Q$. 
The root lattice $L_Q$ associated with $A_Q$ is a free 
abelian group with generators $\alp_1,\dots,\alp_n$, 
called simple roots:
\begin{equation*}
L_Q := \bigoplus_{i=1}^n \Z \alp_i. 
\end{equation*}  
We also define a symmetric bilinear form $(\,,\,) \colon L_Q \times L_Q \to \Z$ by 
$(\alp_i,\alp_j):=(A_Q)_{ij}$. 
It is known that $\chi(S_i,S_j)=(A_Q)_{ij}$, therefore we 
have an isomorphism between lattices
\begin{equation*}
(K(\D_Q),\chi) \cong (L_Q,(\,,\,)), \quad [S_i] \mapsto \alp_i.
\end{equation*}
Reflections 
$r_1,\dots,r_n \colon L_Q \to L_Q$ 
with respect to simple roots $\alp_1,\dots,\alp_n$ 
generates the Weyl group $W:=\left<r_1,\dots,r_n\right>$.  

Let 
$\Delta^{\mathrm{re}}_+ \subset L_Q$ ($\Delta^{\mathrm{im}}_+ \subset L_Q$) be the set of 
positive real (imaginary) roots. 
The imaginary cone $I \subset L_Q \otimes_{\Z} \R$ is 
a closure of a convex hull of 
$\Delta^{\mathrm{im}}_+ \cup \{0\} \subset L_Q \otimes_{\Z} \R$ in the Euclidean topology 
as a finite dimensional vector space. 
Let $V := \Hom_{\Z}(L_Q,\C)$ and introduce the subset $X \subset V$ by
\begin{equation*}
X := V \bs \bigcup_{\lam \in I \bs \{0\}} H_{\lam}
\end{equation*}
where $H_{\lam} := \{\,Z \in V \, \vert \, Z(\lam) = 0\,\}$. 
Further, define the regular subset $\Xreg \subset V$, on which the Weyl group $W$ acts freely, by 
\begin{equation*}
\Xreg := X \bs \bigcup_{\alp \in \Delta_{+}^{\mathrm{re}}} H_{\alp}.
\end{equation*}

A stability condition on $\D_Q$ \cite{Br1} consists of a pair of 
$(Z,\A)$; $\A$ is a full abelian subcategory 
$\A \subset \D_Q$, called 
the heart of a bounded t-structure, and $Z$ is 
a group homomorphism
\begin{equation*}
Z \colon K(\D_Q) \lto \C,
\end{equation*}
called a central charge, which satisfies the condition
\begin{equation*}
Z(E) \in \{\, r e^{i \pi \phi}\ \in \mathbb{C} \,
 \vert \, r \in \mathbb{R}_{>0}, \phi \in (0,1] \,\}
\end{equation*}
for every non-zero object $E \in \A$. 
Let $\Stab(\D_Q)$ be the set of stability conditions on $\D_Q$ with the additional condition, 
called the support property (see Definition \ref{support}). 

The main result in \cite{Br1} implies that $\Stab(\D_Q)$ has 
the structure of a complex manifold and 
there is a local isomorphism map
\begin{equation*}
\pi \colon \Stab(\D_Q) \longrightarrow 
\Hom_{\D_Q}(K(\D_Q),\C), \quad (Z,\A) \mapsto Z
\end{equation*} 
on some open subset of $\Hom_{\D_Q}(K(\D_Q),\C)$. 
Note that under the identification $K(\D_Q) \cong L_Q$, we have 
$\Hom_{\D_Q}(K(\D_Q),\C) \cong V$. 
For the space $\Stab(\D_Q)$, there is a distinguished connected component 
$\Stab^{\circ}(\D_Q) \subset \Stab(\D_Q)$ which contains stability conditions 
with the heart $\A_Q$. 

Seidel and Thomas \cite{ST} defined 
autoequivalences $\Phi_{S_1},\dots,\Phi_{S_n} \in \Aut(\D_Q)$, called spherical twists, 
for spherical objects $S_1,\dots,S_n$, and showed that they satisfy braid relations. 
Denote by $\Br(\D_Q) \subset \Aut(\D_Q)$ the subgroup generated by these spherical twists. 
The action of $\Br(\D_Q)$ on $\Stab(\D_Q)$ preserves 
the distinguished connected component $\Stab^{\circ}(\D_Q)$.

The following theorem is the main result of this paper. 
This generalizes the results by Bridgeland and Thomas \cite{Br3,T}  
for root systems of finite or affine type to arbitrary  
root systems of symmetric Kac-Moody Lie algebras. 
\begin{thm}
\label{main}
There is a covering map 
\begin{equation*}
\underline{\pi} \colon \Stab^{\circ}(\D_Q) \longrightarrow \Xreg \slash W
\end{equation*}
and the subgroup $\Z[2] \times \Br(\D_Q) \subset \Aut(\D_Q)$ acts as the group of 
deck transformations $(\Z[2] \subset \Aut(\D_Q)$ is the subgroup generated by the 
shift functor $[2] \in \Aut(\D_Q)\,)$. 
\end{thm}

By the van der Lek's result \cite{vdL}, the fundamental group of $\Xreg \slash W$ is 
given by
\begin{equation*}
\pi_1(\Xreg \slash W) \cong \Z[\gamma] \times G_W
\end{equation*}
where $G_W=\left<\sigma_1,\dots,\sigma_n\right>$ is the Artin group 
\cite{BS} with generators 
$\sigma_1,\dots,\sigma_n$  
associated with the Weyl group $W=\left<r_1,\dots,r_n\right>$. 
The factor $\Z[\gamma]$ is generated by a loop $\gamma$ around the 
orthogonal hyperplanes of the imaginary cone $I \bs\{0\}$. 
Theorem 1.1 implies that there is a surjective 
group homomorphism 
\begin{equation*}
\widetilde{\rho} \colon \Z[\gamma] \times G_W \to \Z[2] \times \Br(\D_Q).
\end{equation*} 
We can show that $\widetilde{\rho}$ sends the generators $\sigma_1, 
\dots, \sigma_n$ to 
the spherical twists $\Phi_{S_1},\dots,\Phi_{S_n}$ and 
$\gamma$ to the shift functor $[2]$.

The automorphism group $\Aut(\underline{Q})$ of 
the graph $\underline{Q}$ acts on $\A_Q$ by permutating simple modules 
$S_1,\dots,S_n \in \A_Q$ corresponding to vertices of $Q$. 

Let $\Aut^{\circ}(\D_Q) \subset \Aut(\D_Q)$ be the subgroup of 
autoequivalences which preserve the distinguished connected component 
$\Stab^{\circ}(\D_Q)$. Further, write by 
$\Aut_*^{\circ}(\D_Q)=\Aut^{\circ}(\D_Q) \slash \mathrm{Nil}^{\circ}(\D_Q)$ the quotient of 
$\Aut^{\circ}(\D_Q)$ by the subgroup $\mathrm{Nil}^{\circ}(\D_Q)$  consisting of autoequivalences 
which acts trivially on $\Stab^{\circ}(\D_Q)$. 

\begin{cor}
\label{cor_aut}
The group $\Aut_*^{\circ}(\D_Q)$ is given by
\begin{equation*}
\Aut_*^{\circ}(\D_Q) \cong \Z[1] \times ( \Br(\D_Q) \rtimes \Aut(\underline{Q})  ) 
\end{equation*}
where $\Aut(\underline{Q})$ acts on $\Br(\D_Q)$ by permutating the generators 
$\Phi_{S_1},\dots,\Phi_{S_n}$. 
\end{cor}
When $Q$ is of affine ADE type, similar results was shown 
by Corollary 1.5 in \cite{Br3}.

\subsection{Further problems}
Similar to the case of K$3$ surfaces in \cite{Br2} or  
Kleinian singularities in \cite{Br3} (which correspond to 
quivers of finite or affine type), 
we expect the following properties for the space $\Stab(\D_Q)$.

\begin{conjecture}
\begin{itemize}
\label{conjecture}
\item[(1)]
The space $\Stab(\D_Q)$ is connected. Hence $\Stab^{\circ}(\D_Q) =\Stab(\D_Q)$.
\item[(2)]
The space $\Stab^{\circ}(\D_Q)$ is simply connected. In other words, the surjective 
group homomorphism
\begin{equation*}
\rho \colon G_W \longrightarrow \Br(\D_Q)
\end{equation*}
is injective. (Hence $\rho$ is an isomorphism.) 
\end{itemize}
\end{conjecture}

Conjecture \ref{conjecture} (1) was solved for $\hat{A}_1$-quiver in \cite{Ok2}, 
and for $A_n$-quivers and $\hat{A}_n$-quivers in \cite{IUU}. 
Conjecture \ref{conjecture} (2) was solved for $A_n$-quivers in \cite{ST}, 
for ADE-quivers in \cite{BTh}, and for $\hat{A_n}$-quivers in \cite{IUU}. 

Further, the $K(\pi,1)$ conjecture for Artin groups (see \cite{Pa}) together with above 
two conjectures implies that the space $\Stab(\D_Q)$ is contractible. 

Note that if both Conjecture \ref{conjecture} (1) and (2) 
are true, 
then the autoequivalence group of $\D_Q$ is given by
\begin{equation*}
\Aut(\D_Q) \cong \Z[1] \times ( \Br(\D_Q) \rtimes \Aut(\underline{Q})  ). 
\end{equation*}

\subsection*{Acknowledgements}
First of all, I would like to thank my advisor A. Kato for his constant support and 
many valuable discussions. 
I would also like to thank to S. Yanagida for giving many 
valuable comments and pointing out many mistakes. 
I'm also thank to H. Hosaka for many comments and help for 
typesetting the paper. Thanks also to M. Mori, S. Okada, K. Saito 
and Y. Saito 
for helpful discussions and comments. 

\section{Root systems of symmetric Kac-Moody Lie algebras}

\subsection{Root lattices and Weyl groups}
\label{sec:root_lattice}
In this section, we recall basic notions and results for a root system associated with  
a generalized Cartan matrix and describe the Weyl group action on it. 
Further, the sets of real roots and imaginary roots are given. 
We refer to Kac's book and paper \cite{Kac1,Kac2} for more details.

A matrix $A=(a_{ij})_{i,j=1}^n$ is called a generalized Cartan matrix (GCM for short) if $A$ satisfies 
\begin{itemize}
\item[(C1)] $a_{ii}= 2$ for $i=1,\dots , n$,
\item[(C2)] $a_{ij} \in \mathbb{Z}_{\le 0}$ for $i \neq j$,
\item[(C3)] $a_{ij} = 0 \Rightarrow a_{ji} = 0$. 
\end{itemize}
In this paper, we treat only symmetric GCMs, 
so the condition (C3) always true. 

A matrix $A$ is decomposable if $A$ is a block diagonal form
\begin{align*}
A = \left( \begin{array}{lcr}
A_1 & 0 \\
0   & A_2
\end{array} \right)
\end{align*}
up to reordering of indices and indecomposable if otherwise. 

We associate with $A=(a_{ij})_{i,j=1}^n$ a graph $S(A)$, called the 
Dynkin diagram of $A$ as follows. 
The graph $S(A)$ has the set of vertices $\{1,\dots,n\}$, and 
distinct two vertices $i \neq j$ are connected by $\vert a_{ij} \vert$ edges. 
It is clear that $A$ is indecomposable if and only if the graph $S(A)$ is connected. 

For a real column vector $u = (u_1,u_2,\dots)^{T}$, the notation $u >0$ means $u_i >0$ for all $u_i$,  
and $u \ge 0$ means $u_i \ge 0$ for all $u_i$.
The next result gives the classification of indecomposable GCMs. 
\begin{thm}[\cite{Kac1}, Theorem 4.3]
\label{classification}
Let $A$ be an indecomposable GCM. Then one and only one of the following 
three possibilities holds:
\begin{align*}
\text{(Fin) }& \det A \neq 0;\text{there exists }u >0 \text{ such that }
Au >0; Au \ge 0 \Rightarrow u >0 \text{ or } u =0, \\
\text{(Aff) }& \mathrm{rank} A = n-1; \text{there exists }u >0 \text{ such that }
Au =0;Au=0 \Rightarrow u=0, \\
\text{(Ind) }&\text{there exists }u>0 \text{ such that }Au < 0;u \ge 0 \Rightarrow u=0.
\end{align*}
\end{thm}
Referring to cases (Fin), (Aff), or (Ind), we shall say that $A$ is of finite, 
affine, or indefinite type, respectively.

The root lattice associated with $A$ is a free abelian group $L := \oplus_{i=1}^n \mathbb{Z}\alpha_i$ 
with generators $ \alpha_1,\dots, \alpha_n $, called simple roots.  
Denote by $\Pi := \{ \alpha_1,\dots, \alpha_n \}$ the 
set of simple roots. 
Since $A$ is symmetric, we can define a symmetric bilinear form $(\, , \,):L \times L \to \mathbb{Z}$ by 
$(\alpha_i, \alpha_j) := a_{ij}$. 

Define simple reflections $r_i \colon L \to L \,\,\,(i=1,\dots,n)$ by
\begin{equation*}
r_i (\lam) := \lam - (\lam,\alpha_i) \alpha_i, \quad \text{for }\lam \in L. 
\end{equation*}
The group $W := \left< r_1, \dots, r_n \right>$ 
generated by these reflections  
is called the Weyl group 
and satisfies the following relations (see Chapter 3 in \cite{Kac1}):
\begin{equation*}
r_i^2 = 1
\end{equation*}
and for $i \neq j$, 
\begin{align*}
r_i r_j &= r_j r_i \quad \quad \text{if}  \quad  a_{ij} = 0 \\
r_i r_j r_i &= r_j r_i r_j \,\quad \text{if}  \quad  a_{ij} = -1 .
\end{align*}
Note that the symmetric bilinear form $(\, , \,)$ is invariant 
under the  $W$-action;
$(w(\alpha), w(\beta)) = (\alpha, \beta)$ for any $\alpha, \beta \in L$ and $w \in W$. 

For $\alp=\sum_i k_i \alp_i$, the support of $\alp$ is defined to be the full subgraph 
$\supp(\alp) \subset S(A)$ whose vertices are 
$\{\,i \,\vert\, k_i \neq 0 \,\} \subset \{\,1,\dots,n\,\}$. 

Let $L_{+}:=\sum_{i=1}^n \Z_{\ge 0} \alp_i$ and 
$L_{-} := -(L_{+})=\sum_{i=1}^n \Z_{\le 0} \alp_i$. 
We define the set of (positive or negative) 
real roots and imaginary roots by using 
the  $W$-action on $L$. 

\begin{defi} 
\label{root}
\begin{itemize}
\item[(1)]The set of real roots $\Delta^{\mathrm{re}}$ is defined to be the $W$-orbits of simple roots $\Pi$:
\begin{equation*}
\Delta^{\mathrm{re}} := W(\Pi) = \{\,w(\alp_i) \,\vert \, w \in W, i=1,\dots , n \}.
\end{equation*}
The set of positive real roots $\Delta_{+}^{\mathrm{re}}$ (negative real roots $\Delta_{-}^{\mathrm{re}}$) is 
given by
\begin{equation*}
\Delta_{+}^{\mathrm{re}} := \Delta^{\mathrm{re}} \cap L_{+} \quad 
(\Delta_{-}^{\mathrm{re}} := \Delta^{\mathrm{re}} \cap L_{-}\,).
\end{equation*}

\item[(2)]Define the fundamental set of positive imaginary roots $K$ by
\begin{equation*}
K:=\{\,\alp \in L_{+} \bs \{0\} \, \vert \, 
\mathrm{supp}(\alp) \,\,\text{is connected in $S(A)$}, \,
(\alp,\alp_i) \le 0 \,\,\text{for}\,\, i=1, \dots, n     \,\}. 
\end{equation*}
The set of positive imaginary roots $\Delta^{\mathrm{im}}_{+}$ is defined to be the $W$-orbits of $K$:
\[
\Delta^{\mathrm{im}}_{+} := W(K) = \{\,w(\alp) \,\vert \, w \in W, \alp \in K \}.
\]
The set of negative imaginary roots is given by 
$\Delta^{\mathrm{im}}_{-}:= -\Delta^{\mathrm{im}}_{+}$ and 
the set of all imaginary roots is given by 
$\Delta^{\mathrm{im}} := \Delta^{\mathrm{im}}_{+} \cup \Delta^{\mathrm{im}}_{-}$.
\end{itemize}
\end{defi}
\begin{rem}
\label{multiplication}
\begin{itemize}
\item[(1)] For an indecomposable GCM $A$, 
Theorem \ref{classification} implies that the fundamental set $K$ is 
non-empty if and only if $A$ is of affine or indefinite type. 
Hence the set of imaginary roots $\Delta^{\mathrm{im}}$ is also non-empty 
if and only if $A$ is of affine or indefinite type. 
\item[(2)]
Since $K$ is closed under the multiplication of positive integers $\Z_{\ge 1}$ and 
the $W$-action commutes with it, 
the set $\Delta^{\mathrm{im}}_{+} $ is 
also closed under the multiplication of $\Z_{\ge 1}$. 
\end{itemize}
\end{rem}

For a GCM $A=(a_{ij})_{i,j=1}^n$,
take a subset  $J \subset \{1,\dots,n\}$ and  
consider a submatrix $A_J := (a_{ij})_{i,j \in J}$. 
Since $A_J$ is also a GCM, we can define 
the root lattice $L_J$, the Weyl group $W_J$ and 
the set of roots $\Delta_J$ associated with $A_J$. 
By the definitions, we can naturally embed 
$L_J$ and $W_J$ into $L$ and $W$ by  
$L_J := \oplus_{j \in J} \Z \alp_j \subset L$ and  
$W_J:=\left<r_j \,\vert \,j \in J \right> \subset W$. 
For the set of roots $\Delta_J \subset L_J $, 
we can easily see that $\Delta_J = L_J \cap \Delta$.

\subsection{The regular subset $\Xreg$} 
\label{sec_regular}
Let $A$ be a GCM and $L = \oplus_{i=1}^n \Z \alp_i$ be the root lattice associated with $A$.
Throughout this paper, we fix the following notations:
\begin{align*}
V_{\R}^* &:= L \otimes_{\Z} \R = \oplus_{i=1}^n \R \alp_i \\
V^*      &:= L \otimes_{\Z} \C = \oplus_{i=1}^n \C \alp_i = V_{\R}^* \oplus i V_{\R}^* \\
V_{\R} \,&:= \Hom_{\Z}(L,\R)  \\
V \,\,   &:= \Hom_{\Z}(L,\C)  = V_{\R} \oplus i V_{\R}.
\end{align*}

We fix the real structures of complex vector spaces $V$ and $V^*$ as in the above notations. 

The action of the Weyl group $W$ on $L$ is linearly extended to the action on 
$V^* \cong L \otimes_{\Z} \C$ 
(or $V^*_{\R} \otimes_{\Z} \R$) and the contragradient action of $W$ on $V$ (or $V_{\R}$) is defined  
by $\left< w(Z),\lam \right> := \left< Z,w^{-1}(\lam) \right> $ for $Z \in V$ and $\lam \in V^*$.

We fix  some norms for these vector spaces and 
consider the Euclidean topology as finite dimensional vector spaces.

Next we introduce the notion of an imaginary cone which plays central role in this paper. 
\begin{defi}
Assume that the set of positive imaginary roots  $\Delta_{+}^{\mathrm{im}}$ is non-empty. Then, 
the imaginary cone $I$ is defined to be the closure of 
the convex hull of the set $\Delta_{+}^{\mathrm{im}} \cup \{0\}$ in $V^{*}_{\R}$. We  set $I_{0} := I \bs \{0\}$ and 
also call $I_0$ the imaginary cone. 
\end{defi}
For the case $\Delta_{+}^{\mathrm{im}}$ is empty, 
we define that $I$ is empty. 

Note that by Remark \ref{multiplication}, 
the imaginary cone $I$ associated to an indecomposable GCM $A$ 
is non-empty if and only if 
$A$ is of affine or indefinite type. 

\begin{lem}
\label{convex} Assume that the imaginary cone $I$ is non-empty.  
Then, $I$ is a convex cone supported on $\sum_{i=1}^n\,\R_{\ge 0}\alp_i$. 
\end{lem}
{\bf Proof.} The convexity of $I$ is clear by the definition of $I$. 
Since $\Delta_{+}^{\mathrm{im}} \cup \{0\}$ is closed under the multiplication of 
$\Z_{\ge 0}$ (see Remark \ref{multiplication}), 
$I$ is closed under the multiplication of $\R_{\ge 0}$. 
Further $\Delta_{+}^{\mathrm{im}} \cup \{0\}$ is supported on 
$\sum_{i=1}^n\,\Z_{\ge 0}\alp_i$, therefore $I$ 
is supported on $\sum_{i=1}^n\,\R_{\ge 0}\alp_i$. 
\endproof

The following property of the imaginary cone is important 
in the proof of 
Proposition \ref{image}. 
\begin{prop}[\cite{Kac2}, Proposition 1.8 ]
\label{dense}
The set of rays of imaginary roots $\{\, \mathbb{R}_{>0} \alpha \, \vert \, 
\alpha \in \Delta_{+}^{\mathrm{im}} \,\}$ is dense in the imaginary cone $I_0$. 
\end{prop}

For the compatibility of the notations in the later, 
we denote the canonical pairing of $Z \in V$ and $\lam \in V^*$  by 
$Z(\lam) := \left<Z,\lam \right>$.

Now we introduce the open subset $X \subset V $ and 
$\Xreg \subset V$ as follows. 
\begin{defi}
\label{subset}
Let $\lam \in V^*$ and $H_{\lam} := \{\,Z \in V \, \vert \, Z(\lam) = 0\,\}$ denotes 
a complex orthogonal hyperplane with respect to $\lam$.
A subset $X \subset V$ is defined by
\begin{equation*}
X := V \bs \bigcup_{\lam \in I_0} H_{\lam}
\end{equation*}
and a regular subset $X_{\mathrm{reg}} \subset X \subset V$  is defined by
\begin{equation*}
X_{\mathrm{reg}} := X \bs \bigcup_{\alp \in \Delta_{+}^{\mathrm{re}}} H_{\alp}.
\end{equation*}
\end{defi}
Since the $W$-action on $V^*$ preserves real roots $\Delta^{\mathrm{re}}$ 
and the imaginary cone $I_0$, the $W$-action on $X$ and $X_{\mathrm{reg}}$ 
is well-defined.

The following lemma is used to prove that $\Xreg$ is open in $V$.
\begin{lem}[\cite{Kac1}, Lemma 5.8]
\label{lim_ray}
In $V_{\R}^{*} \bs \{0\}$, the limit rays for the set of rays $\{\, \mathbb{R}_{>0} \alpha \, \vert \, 
\alpha \in \Delta_{+}^{\mathrm{re}} \,\}$ lie in $I_0$. 
\end{lem}

Define an imaginary convex disk by 
\begin{equation*}
D:= I \cap \{\,k_1 \alp_1 + \cdots + k_n \alp_n \in V_{\R}^* \,\vert\,
k_1 + \cdots + k_n = 1 \,\}.
\end{equation*}
Lemma \ref{convex} implies that $D$ is a compact convex subset of $V^*_{\R}$. 
  
Note that since $I_0 = \R_{>0} D$, the subset  
$X$ can be written by 
\begin{equation*}
X := V \bs \bigcup_{\lam \in D} H_{\lam}.
\end{equation*}

\begin{lem}
The subsets $X \subset V$ and 
$\Xreg \subset V$ are open in $V$. 
\end{lem}
{\bf Proof.} 
First, we say that 
the subset $X= V \bs \bigcup_{\lam \in D} H_{\lam} 
$ is open. 

Since $D$ is compact, the projectivization $\bigcup_{\lam \in D} 
\mathbb{P} H_{\lam} \subset \mathbb{P}V$ 
is also compact. Hence, the subset 
$\bigcup_{\lam \in D} H_{\lam} $ is closed in $V$. 

Let $Z \in \Xreg$ and take a small open neighborhood 
$Z \in B \subset X$. 
Lemma \ref{lim_ray} implies that the 
accumulated points of the set 
$\bigcup_{\alp \in \Delta^{\mathrm{re}}_+} H_{\alp} $ are 
contained in $\bigcup_{\lam \in D} H_{\lam}$. 
Hence if we take $B$ to be sufficiently small, $B$ 
does not intersect the set 
$\bigcup_{\alp \in \Delta^{\mathrm{re}}_+} H_{\alp} $, and this 
implies that $\Xreg$ is open. 
\endproof

\begin{lem}
\label{lem:image}
Assume that $I_0 \subset V^*$ is non-empty. 
Let $Z \in X$ and consider 
the linear map $Z \colon V^* \to \C$. 
Then the image of the imaginary cone $Z(I_0) \subset \C$ takes the 
following form:
\begin{equation*}
Z(I_0) = \{\,r e^{i \pi \phi} \,\vert\,r >0 \,\,\,\text{and}\,\,\,\phi_1 \le \phi \le \phi_2 \,\}
\end{equation*}
where $\phi_1,\phi_2 \in \R$ (determined up to modulo $2\Z$) with $0 \le \phi_2 - \phi_1 <1 $. 
\end{lem}
{\bf Proof.} Let $D$ be an imaginary convex disk defined above, 
which is compact convex. 
Since $Z \colon V^* \to C$ is $\R$-linear, 
the image $Z(D) \subset \C$ is also 
compact convex and $Z(I_0) = \R_{>0} Z(D)$. 
We also note that $0 \notin Z(I_0)$ since $Z \in X$. 
Thus the set $Z(D)$ is a compact convex subset of
 $\C \bs \{0\}$. 

First we consider the case $Z(D) \cap \R_{>0} = \emptyset$. 
Let  
$\Arg z \in (0,2 \pi)$ for $z \in \C \bs \R_{>0}$ 
the principal argument. 
Then, by the compactness of $Z(D)$, 
we have the maximum and the minimum argument 
of $Z(D)$ by
\begin{align*}
\phi_1 &:= \min \,\left\{\, (1 \slash \pi) \Arg z \in (0,2) \,\vert\,z \in Z(D) \right\} \\
\phi_2 &:= \max \left\{\, (1 \slash \pi) \Arg z \in (0,2) \,\vert\,z \in Z(D) \right\} .
\end{align*}
Clearly $0 \le \phi_2 - \phi_1$, and $\phi_2 - \phi_1 <1$ easily follows from the 
convexity of $Z(D)$ and $0 \notin Z(D)$. 

In the rest case  $Z(D) \cap \R_{>0} \neq \emptyset$, 
the convexity of $Z(D)$ and $0 \notin Z(D)$ implies  
that $Z(D) \cap \R_{<0} = \emptyset$. Therefore the similar argument of the first case works for 
the principal argument $\Arg z \in (-\pi,\pi)$ for $z \in \C \bs \R_{<0}$.  
\endproof

By using the phases $\phi_1 $ and $\phi_2$ which 
appear in 
Lemma \ref{lem:image}, we introduce the phase of the imaginary cone $\phi^I(Z)$ with respect to $Z \in X$, 
up to modulo $2\Z$, by
\begin{align*}
\phi^I(Z):= \frac{\phi_2 + \phi_1}{2}. 
\end{align*}

In the case $I_0$ is non-empty, define the 
normalized regular subset $\Xreg^N$ by 
\begin{equation*}
\Xreg^N := \{\,Z \in \Xreg \,\vert\,\phi^I(Z) = 1 \slash 2 \,\}.
\end{equation*}

By the free action of  $S^1:=\{z \in \C \,\vert\, \abv{z}=1\}$ on $\Xreg$, we have the decomposition  
\begin{align*}
\Xreg \cong S^1 \times \Xreg^N. 
\end{align*}  

Note that the $W$-action on $\Xreg$ preserves the subset $\Xreg^N$ since the $W$-action  
commutes with the $S^1$-action on $\Xreg$. 

\subsection{Fundamental domain}
\label{sec_domain}
In this section, we determine the fundamental domain 
for the action of the Weyl group $W$ on  
the regular subset $X_{\mathrm{reg}}$. To do it, 
we recall some  basic properties of the Tits cones. 
Throughout this section, we assume that the GCM $A$ is indecomposable. 

Recall from Section \ref{sec_regular} that $V$ has the real structure 
$V_{\R} \subset V$ and hence $V=V_{\R} \oplus i V_{\R}$. 
Corresponding to the real structure, $Z \in V$ 
is written by $Z=Z_R + i Z_I$ where $Z_R, Z_I \in V_{\R}$. 

We set $X_{\R} := X \cap V_{\R}$ and  $X_{\R,\mathrm{reg}} 
:= \Xreg \cap V_{\R}$. 
If $A$ is of finite type, then $I=\emptyset$, therefore $X =V$ and $X_{\R} = V_{\R}$. 
If $A$ is of affine or indefinite type, then $X_{\R}$ decomposes 
into two connected components
\begin{align*}
&X_{\R}^+ = \{\,Z_R \in X_{\R}\,\vert
\,Z_{R}(\lam)>0 \quad \text{for all} \quad \lam \in I_0  \} \\
&X_{\R}^- = \{\,Z_R \in X_{\R}\,\vert
\,Z_{R}(\lam)<0 \quad \text{for all} \quad \lam \in I_0  \}.
\end{align*} 
These subsets are known as the Tits cones. 

\begin{defi}[\cite{Kac1}, Section 3.12]
Define the subset $C_{\R} \subset V_{\R}$, called the 
Weyl chamber, by $C_{\R} := \{\,Z_{R} \in V_{\R}\,\vert\,
Z_{R}(\alp_i)>0 \,\,\,\text{for}\,\,\,i=1,\dots,n \}
 (\cong \R^n_{>0})$ and 
let $\overline{C_{\R}} =
\{\,Z_{R} \in V_{\R} \,\vert\,
Z_{R}(\alp_i) \ge 0 \,\,\,\text{for}\,\,\,i=1,\dots,n \} (\cong \R^n_{\ge 0})$ be 
the closure of $C_\R$ in $V_{\R}$. 

The Tits cone $T_{\R}$ is defined by
\begin{equation*}
T_{\R} := \bigcup_{w \in W}\,w(\overline{C_\R}),
\end{equation*}
and the regular subset of $T_{\R}$ is defined by 
\begin{equation*}
T_{\R,\mathrm{reg}} := \bigcup_{w \in W}\,w(C_{\R})
=T_{\R} \bs \bigcup_{\alp \in \Delta_{+}^{\mathrm{re}}} H_{\alp}.
\end{equation*}
\end{defi}

\begin{lem}[\cite{Kac1}, Proposition 3.12 and Section 5.8]
\label{lem:Tits}
The following equality holds:
\begin{align*}
T_{\R} = 
\begin{cases}
V_{\R} \quad &\text{if $A$ is of finite type} \\
X_{\R}^+ \cup \{0\} \quad &\text{if $A$ is of 
affine or indefinite type.}
\end{cases}
\end{align*}
Further,  the Tits cone $T_{\R}$ is a convex cone. 
\end{lem}

\begin{rem}[\cite{Kac1}, Proposition 3.12]
\label{Tits_free}
It is known that the closure of the Weyl chamber $\overline{C_{\R}}$ is a fundamental domain for  
the $W$-action on $T_{\R}$, 
and the $W$-action  on $T_{\R,\mathrm{reg}}$ 
is free and properly discontinuous. 
\end{rem}

The next lemma shall be used in the proof of Proposition \ref{fundamental_domain}.
\begin{lem}
\label{negative}
Assume that a GCM $A$ is of finite type. Then for any 
$Z \in T_{\R,\mathrm{reg}}$, there is 
a unique element $w \in W$ such that $w \cdot Z$ satisfies 
$(w \cdot Z)(\alp_i) < 0$ for all $i = 1,\dots, n$. 
\end{lem}
{\bf Proof.} In the case that $A$ is of finite type, 
Lemma \ref{lem:Tits} implies that 
\begin{align*}
T_{\R,\mathrm{reg}} =V_{\R} \bs \bigcup_{\alp \in \Delta_{+}^{\mathrm{re}}} H_{\alp}.
\end{align*}
Therefore the chamber $- C_{\R} =\{\,Z_R \in V_{R}\,\vert\,
Z_{\R}(\alp_i)<0 \,\,\,\text{for}\,\,\,i=1,\dots,n \}$ is contained in $T_{\R,\mathrm{reg}}$, and is the fundamental domain 
for the action of $W$ on $T_{\R,\mathrm{reg}}$. 
\endproof

In the rest of this section, we assume that 
a GCM $A$ is of affine or 
indefinite type. 
Let $\H:=\{\,z \in \C \,|\,\Im z > 0\,\}$ be the upper half plane and 
define the semi-closed upper half plane by 
\begin{equation*}
H := \H \cup \R_{< 0} = \{\,r e^{i \pi \phi} \in \C \,\vert\, r >0,\,\phi \in (0,1]\,\}.
\end{equation*}
Define the normalized complexified Weyl chamber by 
$C^N := \{\,Z \in \Xreg^N \,\vert \,\, Z(\alp_i) \in H \quad \text{for} \,\,\, i = 1,\dots,n \,\}$.

The following is the main result of this section. 
\begin{prop}
\label{fundamental_domain}
For any $Z \in \Xreg^N$, there is an element $w \in W$ such that 
$w \cdot Z$ lies in $C^N \subset \Xreg^N$.
\end{prop}
{\bf Proof.} Let $Z \in X_{\mathrm{reg}}$. 
Using the real structure $V = V_{\R} \oplus i V_{\R}$, we write
$Z = Z_R + i Z_I $ where $Z_R, Z_I \in V_{\R}$. 
Since $\phi^I (Z) = 1 \slash 2$, the image $Z(I_0)$ is contained 
in $\H$ and this implies that $Z_I(\lam) > 0 $ for all $\lam \in I_0$. 
Hence we have $Z_I \in X_{\R}^+ = T_{\R} \bs \{0\}$.  

Since  $\overline{C_\R}$ is the fundamental domain for the action of 
$W$ on $T_{\R}$, there is some $w^{\prime} \in W$ such that 
$w^{\prime} \cdot Z_I \in \overline{C_\R} $. 

Hence $w^{\prime} \cdot Z = w^{\prime} \cdot Z_R + i w^{\prime} \cdot Z_I$ lies in $V_{\R}+ i \overline{C_\R}$.

We set $Z^{\prime}:=w \cdot Z$, and 
define the subset $J \subset \{1,\dots,n\}$ by
\begin{align*}
J := \{\, j \,\vert\, Z^{\prime}(\alp_j) \in \R \, \}. 
\end{align*}
For $J$, we consider the submatrix $A_J$ as 
in Section \ref{sec:root_lattice}, 
and decompose $A_J$ into indecomposable 
GCMs $A_{J_1},\dots, A_{J_l}$. 

Then, we can say that these indecomposable GCMs $A_{J_1},\dots, A_{J_l}$ are 
of finite type by the following reason. 
If $A_{J_m}\,(m = 1,\dots,l)$ is of affine or indefinite type, 
the image of the imaginary cone $I_0^{J_m} \subset I_0$ 
is also contained in upper half plane; 
$Z^{\prime}(I_0^{J_m}) \subset\H$. 
But by the definition of $J$, $Z^{\prime}(\alp_j) \in \R$ for 
all $j \in J$ and this implies that $Z^{\prime}(I_0^{J_m}) \subset \R$. 
This gives the contradiction, therefore $A_{J_m}$ is of finite type. 

For $Z \colon L \to \C$, consider the restriction 
\begin{equation*}
Z^{\prime} \vert_{L_{J_m}} \colon L_{J_m} \longrightarrow \C
\end{equation*}
where $L_{J_m} \subset L $ is the root lattice associated with  $A_{J_m}$. 

Since $Z^{\prime}(\alp_j) \in \R$ for $j \in J_m $, 
we can regard  $Z^{\prime} \vert_{L_{J_m}}$ as an element of 
the regular subset of the Tits cone $T_{\R,\mathrm{reg}}^{J_m}$ 
associated with $A_{J_m}$.
By Lemma \ref{negative}, we can take an element $w_{J_m} \in W_{J_m}$ such that 
$(w_{J_m} \cdot Z^{\prime})(\alp_j) \in \R_{<0}$ for all $j \in J_m$.

Collect such elements $w_{J_1}, w_{J_2}, \dots,w_{J_l}$ and put  
$w :=w_{J_1}  w_{J_2} \cdots w_{J_l}  w^{\prime}$ ($w$ does not 
depend on the order of 
$w_{J_1}, w_{J_2}, \dots,w_{J_l}$ since these elements commute). 
Then $w$ is the desired element. 
\endproof

\begin{cor}
\label{fundamental_domain2}
For any $Z \in \Xreg$, there are elements $w \in W$ and $k \in \C^{*}$ such that 
$w \cdot k \cdot Z$ lies in $C^N$. 
\end{cor}
{\bf Proof.} Recall from Section \ref{sec_regular} that $\Xreg \cong S^1 \times \Xreg^N$. Hence, for $Z \in \Xreg$, we can take  
$k \in \C^{*}$ such that $Z^{\prime} :=  k \cdot Z$ lies in 
$\Xreg^N$. 
Then, the result follows from Proposition 
\ref{fundamental_domain}. 
\endproof

\begin{prop}
\label{free}
The $W$-action on $\Xreg$ is free and properly discontinuous. 
Further, the fundamental domain for this action is given by 
$S^1 \times C^N \subset S^1 \times \Xreg^N \cong \Xreg$. 
\end{prop}
{\bf Proof.} By Corollary \ref{fundamental_domain2}, 
it is sufficient to consider the element $Z \in C^N$. 
By rotating $Z$, we assume 
that $Z(\alp_i) \in \H $ for all $i = 1,\dots,n $. 
Note that the above condition implies that 
$Z_I \in  C_{\R} \subset T_{\R,\mathrm{reg}}$ 
where  $Z = Z_R + i Z_I$ and $Z_R, Z_I \in V_{\R}$. 
Then, the result follows from the fact that 
the $W$-action on $T_{\R,\mathrm{reg}}$ is 
free and properly discontinuous (see Remark \ref{Tits_free}).  The second 
part immediately follows from Proposition \ref{fundamental_domain}. 
\endproof

\subsection{Walls and chambers in $\Xreg^N$}
\label{wall_chamber}
Here, we introduce the walls in $\Xreg^N$, which are called 
the walls of second kind in \cite{KS}. 
This structure shall be used in 
Section \ref{sec_group} to study the action of the braid group on the space of stability conditions. 

Let $\overline{C^N}$ be a closure of $C^N$ in $\Xreg^N$. 
For $i=1,\dots,n$, 
we define the walls $W_{i,\pm} \subset \overline{C^N}$ by
\begin{align*}
W_{i,+}&:=\{\,Z \in \Xreg^N \,\vert\,Z(\alp_i) \in \R_{>0}, Z(\alp_j) \in \H \text{ for } j \neq i\,\} \\
W_{i,-}&:=\{\,Z \in \Xreg^N \,\vert\,Z(\alp_i) \in \R_{<0}, Z(\alp_j) \in \H \text{ for } j \neq i\,\}.
\end{align*}
Note that $W_{i,-} \subset C^N$, but $W_{i,+} \cap C^N = \emptyset$. 
However $r_i(W_{i,\pm}) = W_{i,\mp} $ and hence $W_{i,+ } \subset r_i(C^N)$. 

By using the $W$-action on $\Xreg^N$, the set of walls 
is defined by
\begin{equation*}
\{\,w(W_{i,\pm})\,\vert\,w \in W, \, i=1,\dots,n\,\}.
\end{equation*}

\begin{lem}
\label{WC}
For any $Z \in W_{i,\pm}$, there is a neighborhood $Z \in U \subset \Xreg^N$ such that 
\begin{equation*}
U \subset C^N \cup r_i(C^N). 
\end{equation*}
\end{lem}
{\bf Proof.} Since $r_i(C^N \cup r_i(C^N)) =C^N \cup r_i(C^N)$ 
and $r_i(W_{i,\pm}) = W_{i,\mp} $, we only need to consider the case that 
$Z \in W_{i,-}$.

Let $Z \in W_{i,-}$. Take an open disk $D_i \subset \C$ centered at $Z(\alp_i) \in \R_{<0}$ such that 
$0 \not\in D_i$ and divide $D_i$ into two pieces 
$
D_{i,+} := \{\,z \in D_i \,\vert\, z \in H\,\}$ and $ D_{i,-} := \{\,z \in D_i \,\vert\, z \in -\H\,\}$. 
Define a neighborhood of $Z$ by $U:=\Xreg^N \cap (D_1 \times \cdots \times D_n) $ where $D_j$ ($j \neq i$) is  
an open disk centered at $Z(\alp_j) \in \H$ which is 
sufficiently small to satisfy $D_j \subset \H$. 
Then, it is easy to check that  
$U_{\pm} :=\Xreg^N \cap (D_1 \times \cdots \times D_{i,\pm} \times \cdots \times D_n)$ 
satisfy $U_+ \subset C^N$ and $U_{-} \subset r_i(C^N)$. 
Hence $U = U_+ \cup U_-$ is the desired open neighborhood. 
\endproof

\subsection{The fundamental group of $X_{\mathrm{reg}}$}
\label{sec_fundamental}
In this section, we give a fundamental group of 
$X_{\mathrm{reg}} \slash W$ by 
using the result of the van der Lek \cite{vdL}. 
This is described in terms of the Artin group associated 
with a Coxter system of the Weyl group $W$ derived from $A$ (\cite{BS}). 
 
\begin{defi}[\cite{BS}]
An Artin group $G_W$ associated with the Weyl group $W$ (derived from $A$) is defined to be the group generated by 
generators $\sigma_1,\dots,\sigma_n$ with the following relations: 
\begin{align*}
\sigma_i \sigma_j &= \sigma_j \sigma_i \quad \quad \text{if}  \quad  a_{ij} = 0 \\
\sigma_i \sigma_j \sigma_i &= \sigma_j \sigma_i  \sigma_j \quad \text{if}  \quad  a_{ij} = -1 .
\end{align*}
\end{defi}

Recall the decomposition $\Xreg \cong S^1 \times \Xreg^N$ in Section \ref{sec_regular}. 
The Weyl group $W$ acts trivially on the first factor $S^1$, hence 
\begin{equation*}
\Xreg \slash W \cong S^1 \times \left(\Xreg^N \slash W \right). 
\end{equation*}
Take a point $*$ in the interior of the chamber $C^N$. 
Let $[*]$ be the class of $*$ in $\Xreg \slash W$.  

\begin{thm}[\cite{vdL}, see also \cite{Pa}, Corollary 3.11]
Assume that a GCM $A$ is of affine or indefinite type. 
Then, the fundamental group for 
$X_{\mathrm{reg}} \slash W$ is given by
\begin{align*}
\pi_1(X_{\mathrm{reg}} \slash W, [*]) \cong 
\Z[\gamma] \times G_W.
\end{align*}
The generator $\gamma$ of the first factor $\Z[\gamma]$ is given by the $S^1$-orbit of $[*]$. 
The generator $\sigma_i$ of the second factor $G_W$ is given by 
the path connecting to $*$ and $r_i(*)$ passing the wall $W_{i,\pm}$ in $\Xreg^N$ just once, 
which is a loop in $\Xreg^N \slash W$. 

\end{thm}
{\bf Proof.} Since
$\Xreg \slash W \cong S^1 \times \left( \Xreg^N \slash W \right)$
and $\pi(S^1) \cong \Z$, it is sufficient to prove that
$\pi_1 (\Xreg^N \slash W) \cong G_W$.
Define the regular subset of the complexified Tits cone by 
\begin{equation*}
T_{\mathrm{reg}}:=\{\,Z \in \Xreg\,\vert\, \Im Z \in T_{\R}\}. 
\end{equation*}
The van der Lek's result in \cite{vdL} implies that $\pi_1(T_{\mathrm{reg}}/W) \cong G_W$. 
Therefore we show that $\Xreg^N$ is homotopic to $T_{\mathrm{reg}}$. 

For $Z \in \Xreg$, $\Im Z \in T_{\R} \bs \{0\} = X_{\R}^+$ is equivalent to 
$Z(I_0) \subset \H$ (see Section \ref{sec_domain}), so we have 
\begin{align*}
T_{\mathrm{reg}}=\left\{\,Z \in \Xreg\,\vert\, 
Z(I_0) \subset \H \,\right\}.
\end{align*}
Construct a deformation retract $h_t \colon T_{\mathrm{reg}} \to T_{\mathrm{reg}}\,$ by 
$h_t(Z) := Z \cdot e^{i \pi t(1/2 - \phi^I(Z))} $ where $0 \le t \le 1 $. 
Then, it is easy to  check $h_1(T_{\mathrm{reg}}) = X^N_{\mathrm{reg}}$ 
and $h_1 = \id $ on $X^N_{\mathrm{reg}} \subset T_{\mathrm{reg}}$.  Hence, it gives a homotopy equivalence 
$T_{\mathrm{reg}} \thicksim X^N_{\mathrm{reg}}$. 
\endproof

\section{Derived categories of preprojective algebras}
\label{sec_derived}
In the following sections, abelian categories and triangulated categories 
are considered to be $\C$-linear categories.

\subsection{Preprojective algebras of quivers}
Let $Q$ be a finite connected quiver without loops. We denote by $Q_0$ its set of vertices 
and $Q_1$ its set of arrows. An opposite quiver $Q^{\op}$ is obtained 
by reversing the orientation of arrows of $Q$. 
For an arrow $a \colon i \to j \in Q_1$, 
we denote the opposite arrow by $a^{*} \colon j \to i \in Q_1^{\op}$. 

A double quiver $\overline{Q}$ is defined by adding all opposite arrows to $Q$, 
so $\overline{Q}_1 = Q_1 \cup Q_1^{\op} $. 
For a quiver $Q$, define an adjacent matrix $(q_{ij})$ of $Q$ by
\begin{equation*}
q_{ij} := \vert \{\,\text{arrows from $i$ to $j$}\,\} \vert.
\end{equation*}

A GCM $A_Q$ associated with $Q$ is defined by 
\begin{equation*}
(A_Q)_{ij} := 2\delta_{ij} - (q_{ij} + q_{ji}) .
\end{equation*}
For a connected quiver $Q$, we say that $Q$ is 
of finite, affine or indefinite type if 
the corresponding GCM $A_Q$ is of finite, affine or indefinite type respectively (see Section \ref{sec:root_lattice}). 

It is known that $A_Q$ is of finite type 
if and only if $Q$ is of ADE type 
and $A_Q$ is of affine type if and only if  $Q$ is of affine ADE type. 

Let us denote by $\C Q$  a path algebra of $Q$ over $\C$. 
We put a gradation on $\C Q$ by using the length of paths. 
\begin{defi}
The preprojective algebra $\Pi (Q)$ associated to $Q$ is defined by 
\begin{equation*}
\Pi (Q) := \C \overline{Q} \slash (\rho) 
\end{equation*}
where $(\rho)$ is an ideal of $\C \overline{Q}$ generated by the element
\begin{equation*}
\rho := \sum_{a \in Q_1} (a a^* - a^* a ).
\end{equation*}
\end{defi}
Since $\rho$ is a homogeneous element in $\C \overline{Q}$, the preprojective algebra 
$\Pi (Q)$ is also a graded algebra by the length of paths. 

Let $\A_Q := \Mod\Pi(Q)$ be an abelian category of finite dimensional nilpotent right $\Pi(Q)$-modules and 
$K(\A_Q)$ be its Grothendieck group.
Let $S_1,\dots,S_n$ 
be simple modules corresponding to the vertices $Q_0=\{1,\dots,n\}$. 
Then the Grothendieck group $K(\A_Q)$ is a free abelian 
group
\begin{equation*}
K(\A_Q) \cong \oplus_{i=1}^n\,\Z\,[S_i]
\end{equation*}
generated by $[S_1],\dots,[S_n]$ which are classes of 
simple modules in $K(\A_Q)$. 

\subsection{Derived categories of preprojective algebras} 
Let $Q$ be a finite connected quiver. In this section, we 
consider the bounded derived category of finite dimensional 
nilpotent modules of the preprojective algebra $\Pi(Q)$. 
\begin{defi}
A triangulated category $\D$ is called $N$-Calabi-Yau $(\CY_N)$ 
if 
for any objects $E,F \in \D$ there is a natural isomorphism
\begin{equation*}
\nu_{E,F}\colon \Hom_{\D}(E,F) \xrightarrow{\sim} \Hom_{\D}(F,E[N])^*
\end{equation*}
where $E[N]$ is the $N$-th shift of an object $E$ and $*$ means the dual complex vector space. 
In other words, the $N$-th shift functor $[N]$ is a Serre functor.
\end{defi}
For $E,F \in \D$, we write $\Hom_{D}^i(E,F) := \Hom_{D}(E,F[i])$. 

\begin{prop}[\cite{Kel}, Section $4$]
Let $D^b(A_Q)$ be the bounded derived category of $\A_Q$. 
If $Q$ is not of finite type, 
then $D^{b}(\A_Q)$ is a $\CY_2$ triangulated category.
\end{prop}
For simplicity, we write $\D_Q := \DQ$.
Since $\D_Q$ is bounded, the Grothendieck group $K(\D_Q)$ of the derived category $\D_Q$ 
is isomorphic to  $K(\A_Q)$ which is defined 
in the last section. 
From now on, we use the notation $K(\D_Q)$ instead of $K(\A_Q)$. 

A bilinear form $\chi \colon \KQ \times \KQ \to \Z$, 
called the Euler form, 
 is defined by 
\begin{equation*}
\chi (E,F) := \sum_{i \in \Z} (-1)^i \dim_{\C} \Hom _{\D_Q}^i(E,F).
\end{equation*}
By the $\CY_2$ property, the Euler form $\chi$ is  
symmetric. 

For the classes of simple modules $[S_i], [S_j] \in \KQ$, 
it is known that the Euler form 
is given by $\chi (S_i,S_j) = a_{ij}$ where $a_{ij}$ is an entry of 
the GCM $A_Q=(a_{ij})$. 

Hence we can identify the lattice $(\KQ,\chi)$ with the root lattice $(L_Q,(\,,\,))$ associated with $A_Q$, through 
the map $[S_i] \mapsto \alp_i$. 

\subsection{Seidel-Thomas braid groups}
We define some autoequivalences of $\D_Q$ which play important role in this paper, called spherical twists, 
introduced by P. Seidel and R.Thomas in \cite{ST}.

An object $S \in \D_Q$ is called $2$-spherical if 
\begin{align*}
\Hom_{\D_Q}^i(S,S) = 
\begin{cases}
\C \quad \text{if} \quad i=0,2  \\
\,0 \quad \text{otherwise} . 
\end{cases}
\end{align*}

\begin{prop}[\cite{ST}, Proposition 2.10]
For a spherical object $S \in \D_Q$, there is an autoequivalence $\Phi_S \in \Aut(\D_Q)$ 
such that there is an exact triangle 
\begin{equation*}
\Hom_{\D_Q}^{\bullet}(S,E) \otimes S \lto E \lto \Phi_S(E)
\end{equation*}
for any object $E \in \D_Q$. The inverse functor $\Phi_S^{-1} \in \Aut(\D_Q)$ is given by
\begin{equation*}
\Phi_S^{-1}(E) \lto E \lto S \otimes \Hom_{\D_Q}^{\bullet}(E,S)^* .
\end{equation*}
\end{prop}
Since simple modules $S_1,\dots,S_n$ of $\A_Q$ 
are $2$-spherical in $\D_Q$, 
they define spherical twists $\Phi_{S_i} \in \Aut(\D_Q)$. 
The Seidel-Thomas braid group $\Br(\D_Q)$ is defined 
to be the subgroup of $\Aut(\D_Q)$ generated by 
spherical twists $\Phi_{S_1},\dots,\Phi_{S_n}$:
\begin{equation*}
\Br(\D_Q) := \left< \Phi_1,\dots,\Phi_n \right>.
\end{equation*}

\begin{prop}[\cite{ST}, Theorem 1.2]
For the group $\Br(\D_Q)$, the following relations hold :
\begin{align*}
\Phi_{S_i}\Phi_{S_j} &= \Phi_{S_j}\Phi_{S_i} \quad \quad \quad \text{if}  \quad  \chi(S_i,S_j) = 0 \\
\Phi_{S_i}\Phi_{S_j}\Phi_{S_i} &= \Phi_{S_j}\Phi_{S_i}\Phi_{S_j} \,\,\quad \text{if}  \quad  \chi(S_i,S_j) = -1.
\end{align*}
\end{prop}

\begin{cor}
There is a surjective group homomorphism 
\begin{equation*}
\rho \colon G_W \to \Br(\D_Q) 
\end{equation*}
defined by $\sigma_i \mapsto \Phi_{S_i}$. 
\end{cor}

Note that at the Grothendieck group level, a spherical twist $\Phi_S$ induces 
a reflection  $[\Phi_S] \colon \KQ \to \KQ$ given by
\begin{equation*}
[\Phi_S]([E]) = [E] - \chi(S,E)[S] 
\end{equation*}
and the inverse functor satisfies $[\Phi_S^{-1}] = [\Phi_S]$. 
In particular, under the identification $(\KQ,\chi) \cong (L_Q,(\, ,\,))$, 
the group $\Br(\D_Q)$ is reduced to the Weyl group $W$ through the map 
$\Phi_{S_i} \mapsto r_i$.

Recall from Section \ref{sec_fundamental} that the fundamental group of $\Xreg \slash W$ is isomorphic to 
$\Z[\gamma] \times G_W$.
We can extend the above group homomorphism $\rho$ to the following group homomorphism
\begin{equation*}
\widetilde{\rho} \colon \Z[\gamma] \times G_W  \lto \Z[2] \times \Br(\D_Q)
\end{equation*}   
which sends $[\gamma]$ to the shift functor $[2] \in \Aut(\D_Q)$. 

\section{Bridgeland stability conditions}
In the following sections, we always assume that 
the Grothendieck group $K(\D)$ of a triangulated category $\D$ 
is free of finite rank ($K(\D) \cong \Z^n$ for some $n$).  
\label{sec:stab}
\subsection{The spaces of stability conditions}

In this section, we recall the notion of 
stability conditions on triangulated categories 
introduced by T. Bridgeland in \cite{Br1} and collect 
some basic results for the space of stability conditions following \cite{Br1,Br2,BrSm}. 

\begin{defi}
Let $\D$ be a triangulated category and $K(\D)$ be its $K$-group. 
A stability condition $\sigma = (Z, \sli)$ on $\D$ consists of 
a group homomorphism $Z \colon K(\D) \to \C$ called central charge and 
a family of full additive subcategories $\sli (\phi) \subset \D$ for $\phi \in \R$ 
satisfying the following conditions:
\begin{itemize}
\item[(a)]
if  $0 \neq E \in \sli(\phi)$, then 
$Z(E) = m(E) \exp(i \pi \phi)$ for some $m(E) \in \R_{>0}$, 
\item[(b)]
for all $\phi \in \R$, $\sli(\phi + 1) = \sli(\phi)[1]$, 
\item[(c)]if $\phi_1 > \phi_2$ and $A_i \in \sli(\phi_i)\,(i =1,2)$, 
then $\Hom_{\D}(A_1,A_2) = 0$,
\item[(d)]for $0 \neq E \in \D$, there is a finite sequence of real numbers 
\[
\phi_1 > \phi_2 > \cdots > \phi_m
\]
and a collection of exact triangles
\begin{equation*}
0 =
\xymatrix{ 
 E_0 \ar[rr]   &&  E_1 \ar[dl] \ar[rr] && E_2 \ar[dl] 
 \ar[r] & \dots  \ar[r] & E_{m-1} \ar[rr] && E_m \ar[dl] \\
& A_1 \ar@{-->}[ul] && A_2 \ar@{-->}[ul] &&&& A_m \ar@{-->}[ul] 
}
= E
\end{equation*}
with $A_i \in \sli(\phi_i)$ for all $i$.
\end{itemize}
\end{defi}

It follows from the definition that the subcategories $\sli(\phi) \subset \D$ are abelian categories 
(see Lemma 5.2. in \cite{Br1}). 
The nonzero objects of $\sli(\phi)$ are called semistable of phase $\phi$ in $\sigma$, and simple objects 
in $\sli(\phi)$ are called stable of phase $\phi$ in $\sigma$. 

For a stability condition $\sigma = (Z,\sli)$, 
we introduce the set of semistable classes $\class^{\mathrm{ss}}(\sigma) \subset K(\D)$ by
\begin{equation*}
\class^{\mathrm{ss}}(\sigma) :=\{\,\alp \in 
K(\D)\,\vert\,\text{there exists a semistable object } E \in \D \text{ in $\sigma$ such that } [E] = \alp\,\}.
\end{equation*} 
Similarly the set of stable classes $\class^{\mathrm{s}}(\sigma)$ can be defined.

We always assume our stability conditions satisfy the additional assumption 
called the support property in \cite{KS}. 
\begin{defi}
\label{support}
Let $\norm{\,\cdot\,}$ be some norm on $K(\D) \otimes \R$. A stability condition $\sigma=(Z,\sli)$ 
has a support property if there is a some constant $C >0$ such that 
\begin{equation*}
C \cdot \abv{Z(\alp)} > \norm{\alp}
\end{equation*}
for all $\alp \in \class^{\mathrm{ss}}(\sigma)$. 
\end{defi}

\begin{rem}
For a ray $R = \R_{>0} \alp \subset K(\D) \otimes \R$ (where $\alp \in K(\D) \bs \{0\}$), 
define a function $f$ by 
\begin{equation*}
f(R ):= \frac{\abv{Z(\alp)}}{\norm{\alp}}.
\end{equation*}   
(This doesn't depend on the choice of $\alp$, only depends on the ray.)

Let $\R_{>0} \class^{\mathrm{ss}}(\sigma) 
:=\{\,\R_{>0}\alp\,\vert\,\alp \in \class^{\mathrm{ss}}(\sigma)\,\}$ be the set of rays generated by 
semistable classes of $\sigma$. 
Then, the support property of $\sigma =(Z,\sli)$ is equivalent to that 
there is no sequence of rays $R_i\subset \R_{>0} \class^{\mathrm{ss}}(\sigma)\,(i=1,2,\dots)$
such that 
\begin{equation*}
\lim_{i \to \infty} f(R_i) = 0.
\end{equation*}
\end{rem}
 
Let $\Stab(\D)$ be the set of all stability conditions on $\D$ with the support property. 
In \cite{Br1}, Bridgeland introduced a natural topology on $\Stab(\D)$ induced by the metric 
$d \colon \Stab(\D) \times \Stab(\D) \to [0,\infty]$. 
For more details, we refer to Section 8 in \cite{Br1}.  

In this topology, Bridgeland showed the following crucial theorem. 
\begin{thm}[\cite{Br1}, Theorem 1.2]
\label{local_iso}
The space $\Stab(\D)$ has the structure of a complex manifold and 
the projection map of central charges
\begin{equation*}
\pi \colon \Stab(\D) \longrightarrow \Hom_{\Z}(K(\D),\C)
\end{equation*}
defined by $(Z,\sli) \mapsto Z$ is a local isomorphism onto an open subset of $\Hom_{\Z}(K(\D),\C)$.
\end{thm}

The next lemma implies local injectivity of the above
projection map $\pi$. 
\begin{lem}[\cite{Br1}, Lemma 6.4]
\label{distance}
Let $\sigma =(Z,\sli)$ and $\sigma^{\prime} = (Z,\sli^{\prime})$ be stability conditions 
on $\D$ with the same central charge $Z$. 
Then, $d(\sigma,\sigma^{\prime}) < 1$ implies $\sigma = \sigma^{\prime}$. 
\end{lem}

The following lemma shall be used in the proof 
of  Proposition \ref{ssc} and Proposition \ref{image}. 
\begin{lem}
\label{closed}
Fix a class $\alp \in K(\D)$ and let $U \subset \Stab(\D)$ be an open subset. 
If every stability condition $\sigma \in U$ 
satisfies $\alp \in \class^{\mathrm{ss}}(\sigma)$, 
then a stability condition on the boundary $\sigma^{\prime} \in \rd U$ 
also satisfies $\alp \in \class^{\mathrm{ss}}(\sigma^{\prime})$.
\end{lem}
{\bf Proof.}
This follows from the results for walls and chambers in \cite[Section 9]{Br2} 
or \cite[Section 7.6]{BrSm}.
\endproof

\subsection{Stability conditions on finite length abelian categories}
\label{finite_length}
In \cite{Br1}, Bridgeland gave the alternative description of a stability condition on $\D$ 
as the pair of a bounded t-structure and a central charge on its heart. In this section, by using this description, we construct stability conditions 
on finite length abelian categories with finitely many 
simple objects. 

\begin{defi}
\label{central_charge}
Let $\mathcal{A}$ be an abelian category and let $K(\mathcal{A})$ be its Grothendieck group. 
A central charge on $\mathcal{A}$ is 
a group homomorphism $Z \colon K(\mathcal{A}) \to \mathbb{C}$ such that for any nonzero object
$0 \neq E \in \mathcal{A}$, the complex number $Z(E)$ lies in semi-closed upper half-plane 
$H = \{\, r e^{i \pi \phi}\ \in \mathbb{C} \,
 \vert \, r \in \mathbb{R}_{>0}, \phi \in (0,1] \,\}$.
\end{defi}

The real number 
$\phi(E):= (1 \slash \pi)\,\arg Z(E) \in (0,1]$
for $0 \neq E \in \mathcal{A}$ is called the phase of $E$. 

A nonzero object $0 \neq E \in \mathcal{A}$ is said to be $Z$-(semi)stable if 
for any nonzero proper subobject $0 \neq A \subsetneq E$ satisfies $\phi(A)  < (\le) \phi(E)$.

\begin{prop}[\cite{Br1}, Proposition 5.3]
\label{abel_stability}
Let $\D$ be a triangulated category. To give a stability condition on $\D$ is equivalent 
to giving the heart $\A \subset \D$ 
of a bounded structure on $\D$ and a central charge with the Harder-Narasimhan property on $\A$. 
\end{prop}
For the  heart of a bounded t-structure 
and the Harder-Narasimhan property (HN property), we refer 
to Section 2 and 3 in \cite{Br1}. 

We denote by $\Stab({\mathcal{A}})$ the set of central charges on 
the heart  $\A \subset \D$  
with the HN property and the support property. 

Proposition \ref{abel_stability} implies that there is a natural inclusion 
\begin{equation*}
\Stab(\A) \subset \Stab(\D). 
\end{equation*}  

Let $\A \subset \D$ be a finite heart, which is an abelian category with finitely many simple objects $S_1,\dots,S_n$ and 
generated by means of extensions of these simple objects. 
Then, we have $K(\A) \cong \oplus_{i=1}^n \Z [S_i]$. 
For any point $(z_1,\dots,z_n) \in H^n$, the central charge $Z \colon K(\A) \to \C$ is 
defined by $Z(S_i) :=z_i$. Conversely, for a given central charge 
$Z \colon K(\A) \to \C$, the complex number 
$Z(S_i)$ lies in $H$ for all $i$. Hence $Z$ determines 
a point $(z_1,\dots,z_n) \in H^n$ where $z_i :=Z(S_i)$.
As a result, the set of central charges on $\A$ is 
isomorphic to $H^n$.

\begin{lem}[\cite{Br4}, Lemma 5.2]
Let $Z \colon K(\A) \to \C$ be a central charge given by the above construction. 
Then $Z$ has the HN property. 
In particular, we have 
\begin{equation*}
\Stab(\A) \cong H^n. 
\end{equation*}
\end{lem}

As in Section \ref{sec_derived}, let $\A_Q$ be an abelian category of finite dimensional nilpotent 
$\Pi(Q)$-modules and $\D_Q$ be the bounded derived category
of $\A_Q$. 
Then, there is a unique distinguished connected component 
$\Stab^{\circ}(\D_Q) \subset \Stab(\D_Q)$ which contains 
the subset $\Stab(\A_Q)$.

\subsection{Group actions on $\Stab^{\circ}(\D_Q)$}
\label{sec_group}
Here we consider two group actions on $\Stab^{\circ}(\D_Q)$, which 
are lift of the action of $\C^*$ and $W$ on $\Xreg$. 
Further, we prove the lifted version of Proposition \ref{fundamental_domain} and Corollary \ref{fundamental_domain2}. 

For the space $\Stab(\D)$, 
we introduce two group actions which commute each other. 

First, consider the action of $\Aut(\D)$ on $\Stab(\D)$. 
Let $\Phi \in \Aut(\D)$ be an autoequivalence of $\D$ and 
 $(Z,\sli) \in \Stab(\D)$ be a stability condition on $\D$. Then, the element 
$\Phi \cdot (Z,\sli) =(Z^{\prime},\sli^{\prime})$ is defined by  
\begin{equation*}
Z^{\prime}(E):=Z(\Phi^{-1}(E)), \quad \sli^{\prime}(\phi) := \Phi(\sli(\phi)),
\end{equation*}
where $E \in \D$ and $\phi \in \R$.
 
Secondly, define the $\C$-action on $\Stab(\D)$. 
For $t \in \C$ and $(Z,\sli) \in \Stab(\D)$, 
the element $t \cdot (Z,\sli) =(Z^{\prime},\sli^{\prime})$ is defined by 
\begin{equation*}
Z^{\prime}(E):= e^{-i \pi t} \cdot Z(E), \quad \sli^{\prime}(\phi) :=\sli(\phi + \Re (t))
\end{equation*}
where $E \in \D$ and $\phi \in \R$. 
Clearly, this action is free. 

Note that by the definition, these two actions are 
isometries with respect to the distance $d$ on the space $\Stab(\D)$.  

On $K(\D)$, we can easily see the following remark about semistable classes. 
\begin{rem}
\label{groupss}
By the action of an autoequivalence $\Phi \in \Aut (\D)$, the set of semistable classes 
$\class^{\mathrm{ss}}(\sigma)$ for 
$\sigma \in \Stab(\D)$ changes to
\begin{equation*}
\class^{\mathrm{ss}}(\Phi \cdot \sigma) = [\Phi](\class^{\mathrm{ss}}(\sigma) ).
\end{equation*}
On the other hand, the $\C$-action on $\Stab(\D)$ 
does not change the set of 
semistable classes $\class^{\mathrm{ss}}(\sigma)$. 
\end{rem}

In the rest of this section, we study the action of the Seidel-Thomas 
braid group $\Br(\D_Q)$ on $\Stab^{\circ}(\D_Q)$. 
Recall the identification $K(\D_Q) \cong L_Q$ and consider the 
projection map 
\begin{equation*}
\pi \colon \Stab^{\circ}(\D_Q) \longrightarrow V
\end{equation*}
where $V = \Hom_{\Z}(L_Q,\C)$ (see Section \ref{sec_regular}). 

Then, $\pi $ maps the subset $H^n \cong \Stab(\A_Q) \subset \Stab^{\circ}(\D_Q)$ 
isomorphically onto the subset 
\begin{equation*}
 H^n \cong \{\,Z \in \Xreg \,\vert \,\, Z(\alp_i) \in H \quad \text{for} \,\,\, i = 1,\dots,n \,\} \subset V.
\end{equation*}

Corresponding to the normalized regular subset $\Xreg^N \subset V$, 
we introduce the space of normalized stability conditions $\Stab(\D_Q)^N$ by
\begin{equation*}
\Stab(\D_Q)^N:= \{\,\sigma \in \Stab^{\circ}(\D_Q)\,\vert\,\pi(\sigma ) \in \Xreg^N\,\},  
\end{equation*}
and normalized stability conditions on $\A_Q$ by 
\begin{equation*}
\Stab(\A_Q)^N:=\{\,\sigma \in \Stab(\A_Q)\,\vert\,\pi(\sigma ) \in \Xreg^N\,\}. 
\end{equation*}
Note that the projection  
$\pi \colon \Stab(\D_Q)^N \to \Xreg^N$
maps $\Stab(\A_Q)^N$ isomorphically onto the chamber $C^N \subset \Xreg^N$. 

Recall from Section \ref{wall_chamber} that 
there are walls $W_{i,\pm} \subset \overline{C^N}$ for $i=1,\dots,n$. 
Define the lifted walls $\widetilde{W_{i,\pm}} \subset \overline{\Stab(\A_Q)^N}$ for $i=1,\dots,n $ by 
\begin{align*}
\widetilde{W_{i,+}}&:=\{\,\sigma=(Z,\sli) \in \overline{\Stab(\A_Q)^N}\,\vert\,
Z(S_i) \in \R_{>0}, Z(S_j) \in \H \text{ for } j \neq i\,\} \\
\widetilde{W_{i,-}}&:=\{\,\sigma=(Z,\sli) \in \overline{\Stab(\A_Q)^N}\,\vert\,
Z(S_i) \in \R_{<0}, Z(S_j) \in \H \text{ for } j \neq i\,\}.
\end{align*}
Note that as in Section \ref{wall_chamber}, 
$\widetilde{W_{i,-}} \subset \Stab(\A_Q)^N$ but 
$\widetilde{W_{i,+}} \cap \Stab(\A_Q)^N = \emptyset$. However, 
$\Phi_{S_i}^{-1}(\widetilde{W_{i,-}})=\widetilde{W_{i,+}}$ and 
$\widetilde{W_{i,+}} \subset \Phi_{S_i}^{-1}(\Stab(\A_Q)^N)$.

The following is the lift of Lemma \ref{WC}. 
\begin{lem}
\label{WC2}
Let $\sigma \in \widetilde{W_{i,\pm}} \subset \overline{\Stab(\A_Q)^N}$. Then, there is a 
neighborhood $\sigma \in U \subset \Stab(\D_Q)^N$ such that one of 
the following holds
\begin{itemize}
\item[(1)] $U \subset \Stab(\A_Q)^N \cup \Phi_{S_i}^{-1}(\Stab(\A_Q)^N)$ 
if $\sigma \in \widetilde{W_{i,+}}$, 
\item[(2)] $U \subset \Stab(\A_Q)^N \cup \Phi_{S_i}\:(\Stab(\A_Q)^N)$ 
if $\sigma \in \widetilde{W_{i,-}}$.
\end{itemize}
\end{lem}
{\bf Proof.} Note that in a $\CY_2$ category, 
simple tilted categories $\mu_{S_i}^{\pm}(\A_Q)$ correspond 
to $\Phi_{S_i}^{\pm 1}(\A_Q)$. Then it follows from Lemma 5.5 in \cite{Br4} or Lemma 7.9 in \cite{BrSm}. 
\endproof

\begin{lem}
\label{contain}
The image of the projection map $\pi \colon \Stab^{\circ}(\D_Q) \to V $
contains $\Xreg$.
\end{lem}
{\bf Proof.} Recall from Corollary \ref{fundamental_domain2} that 
the orbit of $C^N \subset V$ under the 
action of $\C^*$ and $W$ coincides with $\Xreg$. 
Since the action of $\C$ and $\Br(\D_Q)$ on 
$\Stab^{\circ}(\D_Q)$ is reduced to 
the action of $\C^*$ and $W$ on the base space $V$, 
the orbit of $\Stab(\A_Q)^N \subset \Stab^{\circ}(\D_Q)$ under the action 
of $\C$ and $\Br(\D_Q)$ is mapped to the subset $\Xreg \subset V$.
\endproof

Let $\Stab^{\circ}(\D_Q)^N$  be the connected component of $\Stab(\D_Q)^N$ 
which contains $\Stab(\A_Q)^N$. 
Now, we lift Proposition \ref{fundamental_domain} and Corollary \ref{fundamental_domain2} via 
the restricted projection map $\pi \colon \pi^{-1}(\Xreg) \to \Xreg$. 

In the proof of the following result, we use the same argument in  
the proof of Proposition 13.2 in \cite{Br2}. 
\begin{prop}
\label{orbit}
For any $\sigma \in \Stab^{\circ}(\D_Q)^N$, there is an autoequivalence $\Phi \in \Br(\D_Q)$ such 
that $\Phi \cdot \sigma$ lies in $\Stab(\A_Q)^N$.  
\end{prop}
{\bf Proof.} Let $\sigma=(Z,\sli) \in \Stab^{\circ}(\D_Q)^N$ and 
take a path $\gamma \colon [0,1] \to \Stab^{\circ}(\D_Q)^N$ 
such that $\gamma (0) \in \Stab(\A_Q)^N$ and $\gamma(1) = \sigma$. 
By Lemma \ref{WC}, 
we can deform $\gamma$ to satisfy that for any $t \in (0,1)$, 
the path 
$\pi(\gamma) \subset \Xreg^N$ passes 
the walls $\{\,w(W_{i,\pm})\,\vert\,w \in W, \, i=1,\dots,n\,\}$ only at 
$t_1,\dots,t_m \in (0,1)$ with
$0<t_1< \cdots <t_m <1$. 

Since $\gamma(\,[0,t_1)\,) \subset \Stab(\A_Q)^N $ and $\pi(\gamma(t_1)) \in W_{i,\pm}$ for some $i$, 
the stability $\gamma(t_1) \in \Stab^{\circ}(\D_Q)^N$ lies in $\widetilde{W_{i,\pm}}$. 
If $\gamma(t_1) \in \widetilde{W_{i_1,+}}$, define $\Phi_1 :=\Phi_{S_{i_1}}$, 
and if $\gamma(t_1) \in \widetilde{W_{i,-}}$, define $\Phi_1 :=\Phi_{S_{i_1}}^{-1}$. 
Then, by Lemma \ref{WC2}, $(\Phi_1 \gamma)(\,(t_1,t_2)\,) \subset \Stab(\A_Q)^N$ and 
$\pi((\Phi_1 \gamma)(t_2)) \in W_{i_2,\pm}$. 
Hence $(\Phi_1  \gamma)(t_2) \in \widetilde{W_{i_2,\pm}}$, and 
we can similarly define  
$\Phi_2 :=\Phi_{S_{i_2}}\Phi_1$ or $\Phi_2 := \Phi_{S_{i_2}}^{-1} \Phi_1$. 

Repeating this process, we get an autoequivalence $\Phi_m \in \Br(\D_Q)$ 
such that $(\Phi_m  \gamma)(\,(t_m,1]\,) \subset \Stab(\A_Q)^N$. 
\endproof

Let $\pi^{-1}(\Xreg)^{\circ}$ be the connected component of $\pi^{-1}(\Xreg)$ 
which contains $\Stab(\A_Q)$. 

\begin{cor}
\label{action}
For $\sigma \in \pi^{-1}(\Xreg)^{\circ} \subset \Stab^{\circ}(\D_Q)$, 
there are elements $\Phi \in \Br(\D_Q)$ and $k \in \C$ such that 
$\Phi \cdot k \cdot \sigma \in  \Stab(\A_Q)^N$. 
\end{cor}
{\bf Proof.} Let $\sigma \in \pi^{-1}(\Xreg)^{\circ}$ and take a path 
$\gamma \colon [0,1] \to \pi^{-1}(\Xreg)^{\circ}$ 
such that $\gamma (0) \in \Stab(\A_Q)$ and $\gamma(1) = \sigma$. 
By the $\C$-action on $\pi^{-1}(\Xreg)^{\circ}$, we can normalize $\gamma$ 
to the path $\gamma^{\prime}=k \cdot \gamma$ which lies in $\Stab^{\circ}(\D_Q)^N$ 
where $k \colon [0,1] \to \C$ and $\gamma^{\prime}(t) = k(t)\cdot \gamma(t)$. 
Then, the result follows from Proposition \ref{orbit}.
\endproof

\section{Proof of main theorem}
\subsection{Indivisible roots and semistable classes}
In this section, we show that the set of indivisible roots 
are contained in the set of semistable classes of 
the stability condition with the central charge in $\Xreg$.  

Set $K(\D_Q)_{\ge 0}:=\sum_{i=1}^n \,\Z_{\ge 0}[S_i]$ 
and $K(\D_Q)_{>0}:=K(\D_Q)_{\ge 0}\bs\{0\}$. 
For  classes $\alp,\beta \in K(\D_Q)$, we write by 
$\alp > \beta$  if $\alp - \beta \in K(\D_Q)_{> 0}$. 

In the following we fix the class $\alp \in K(\D_Q)_{> 0}$. We 
denote by $\Rep(\overline{Q},\alp)$ the 
affine space consisting of representations of the double quiver $\overline{Q}$ with the class $\alp$.  
Let 
$\Rep(\Pi(Q),\alp)^{\nil}$ be the set of nilpotent 
$\Pi(Q)$-modules (representations) with the class 
$\alp$.  Lusztig showed that $\Rep(\Pi(Q),\alp)^{\nil}$ 
is a Lagrangian subvariety of the affine space 
$\Rep(\overline{Q},\alp)$ (see Section 12 in \cite{Lus}).

Let $Z \colon K(\D_Q) \to \C$ be a central charge on $\A_Q$. 
The central charge $Z$ is identified with the King's 
stability condition  $\lam \colon K(\D_Q) \to \R$ \cite{King}
by defining    
\begin{equation*}
\lam(\beta):=-\Im\frac{Z(\beta)}{Z(\alp)}, \quad \beta \in K(\D_Q).
\end{equation*}
Therefore, here we use the notion of a central charge 
instead of the notion of a King's stability condition which 
is used in \cite{CBvdB}. 

Following \cite{CBvdB}, we introduce generic stability conditions. 
\begin{defi} 
A central charge $Z \colon K(\D_Q) \to \C$ 
is said to be generic with respect to $\alp$ 
if $\Im(Z(\beta) \slash Z(\alp)) \neq 0$ for all $0 < \beta < \alp$.  
\end{defi}

Let $\mathfrak{g}(A_Q)$ be the Kac-Moody Lie algebra 
associated with the GCM $A_Q$ (see Chapter 1 in \cite{Kac1}). 
The root multiplicity of a root $\alp \in \Delta$ is defined to be 
the dimension $\dim \mathfrak{g}_{\alp}$ where 
$\mathfrak{g}_{\alp}$ is the root space given by 
the 
decomposition $\mathfrak{g}(A_Q)=
\mathfrak{h}\oplus \oplus_{\alp \in \Delta}\mathfrak{g}_{\alp}$.

A root $\alp \in \Delta$ is called indivisible if there is no $\beta \in \Delta$ satisfying 
$\alp = m \beta $ for $\abv{m} > 1$. 

\begin{prop}[\cite{CBvdB}, Proposition 1.2]
\label{indivisible_root}
Let $\alp \in \Delta_{+}$ be a positive indivisible root and suppose that $Z$ is generic with respect to $\alp$. 
Then, the number of irreducible components of $\Rep(\Pi (Q),\alp)^{\mathrm{nil}}$ which contain a $Z$-stable 
representation is 
equal to the root multiplicity $\dim \mathfrak{g}_{\alp}$. 
\end{prop}

As a result, it turns out that 
if $Z$ is generic with respect to a positive indivisible root $\alp \in \Delta_+$, 
then the moduli space of $Z$-stable nilpotent 
modules with the class $\alp$ is non-empty. 

The following is a key result in the proof of 
Proposition \ref{image}.
\begin{prop}
\label{ssc}
Let $\sigma=(Z,\sli) \in \pi^{-1}(\Xreg)^{\circ}$. 
Then, the set of indivisible roots in $\Delta$ is contained in the set of $\sigma$-semistable classes 
$\class^{\mathrm{ss}}(\sigma)$ :
\begin{equation*}
\{ \,\alp \in \Delta \,\vert \, \alp \text{ is indivisible}\} 
\subset \class^{\mathrm{ss}}(\sigma).
\end{equation*}
\end{prop}
{\bf Proof.} First note that 
the set of all indivisible roots are invariant under the action of 
$W$. Since $\sigma \in \pi^{-1}(\Xreg)^{\circ}$, 
by Remark \ref{groupss} and Corollary \ref{action}, it is sufficient to prove 
that any stability condition $\sigma =(Z,\sli) \in \Stab(\A_Q)$ contains 
all indivisible roots as semistable classes. 
 
Fix a positive indivisible root $\alp \in \Delta_+$ and consider the dense open subset of 
$\Stab(\A_Q)$ consisting of stability conditions which are generic with respect to $\alp$: 
\begin{equation*}
\{\,Z \in \Stab(\A_Q)\,\vert\, \Im ( Z(\beta) \slash  Z(\alp) )
\neq 0\text{ for all $0 < \beta < \alp $}\}.
\end{equation*}

By Proposition \ref{indivisible_root}, any stability condition in this subset contains $\alp$ as a stable class. 
Since this subset is dense in $\Stab(\A_Q)$, by Lemma \ref{closed} we conclude that 
any stability condition in $\Stab(\A_Q)$ contains $\alp$ at least as a semistable class. 
\endproof

\subsection{Projection of central charges}
In this section, we determine the image of central charges via the projection map 
\begin{equation*}
\pi \colon \Stab^{\circ}(\D_Q) \to V
\end{equation*}
by using Proposition \ref{ssc} proved in the last section. 

\begin{lem}
\label{boundary}
Let $\rd \Xreg$ be a boundary of $\Xreg$ and assume that $Z \in \rd \Xreg$. 
Then, there is at least one ray $R \subset I_0 \cup \R_{>0} \Delta^{\mathrm{re}}_+$ such that 
$Z(R) = 0$. 
\end{lem}
{\bf Proof.} It immediately follows from the definition of the open subset 
$\Xreg \subset V$ (see Definition \ref{subset}).
\endproof
 
\begin{prop}
\label{image}
The projection map 
\begin{equation*}
\pi \colon \Stab^{\circ}(\D_Q) \to V
\end{equation*} 
maps $\Stab^{\circ}(\D_Q)$ onto the subset $\Xreg \subset V $.
\end{prop}
{\bf Proof.} One inclusion $\Xreg \subset \pi(\Stab^{\circ}(\D_Q))$ is 
Lemma \ref{contain}. 
Here we prove the other inclusion $\pi(\Stab^{\circ}(\D_Q)) \subset \Xreg$,  
which is equivalent to that $\pi^{-1}(\Xreg)^{\circ}= \Stab^{\circ}(\D_Q)$. 

Since $\Stab^{\circ}(\D_Q)$ is the connected component which contains $\pi^{-1}(\Xreg)^{\circ}$, 
it is sufficient to prove that $\pi^{-1}(\Xreg)^{\circ}$ is open and closed. 
 
First note that since $\Xreg$ is open, the connected component $\pi^{-1}(\Xreg)^{\circ}$ is also open. 
Hence, the closedness of $\pi^{-1}(\Xreg)^{\circ}$ is equivalent to that it has no boundary points. 

Assume that $\pi^{-1}(\Xreg)^{\circ}$ has a boundary point $\sigma=(Z,\sli)$. 
Then, $\sigma = (Z,\sli)$ is projected on $\rd \Xreg$. Hence, 
by Lemma \ref{boundary}, 
there is a ray $R \subset I_0 \cup \R_{>0} \Delta^{\mathrm{re}}_+$ such that $Z(R) = 0$. 

For the ray $ R \subset I_0 \cup \R_{>0} \Delta^{\mathrm{re}}_+$, 
by Proposition \ref{dense}, we can take a sequence of rays $R_i\,(i=1,2,\dots,)$ 
such that $R_i \to R$ (as $i \to \infty$) 
where $R_i = \R_{>0} \alp_i$ and each $\alp_i \in \Delta_+$ is an indivisible 
positive root. 

On the other hand, since $\sigma$ lies in the closure of $\pi^{-1}(\Xreg)^{\circ}$, 
by Lemma \ref{closed} and Proposition \ref{ssc}, 
$\sigma$ contains all indivisible roots as semistable classes.   
In particular, the above rays $R_i\, (i=1,2,\dots)$ are contained in 
$\R_{>0} \class^{\mathrm{ss}}(\sigma)$. 

Since $Z(R)=0$, we have 
\begin{equation*}
\lim_{i \to \infty}f(R_i) =f(R) =0
\end{equation*}
where $f$ is the function defined in Remark \ref{support}.
But this contradicts to the support property of $\sigma$ (see Remark \ref{support}). 
\endproof

\subsection{Covering structures}

\begin{prop}
The action of $\Z[2] \times \Br(\D_Q) \subset \Aut(\D_Q)$ on 
$\Stab^{\circ}(\D_Q)$ is  free and properly discontinuous.
\end{prop}
{\bf Proof.} 
The result is clear for the action of $\Z[2]$. Hence we prove 
them for the action of $\Br(\D_Q)$.

We first prove that the action of $\Br(\D_Q)$ is free. By Corollary \ref{action}, 
it is sufficient to prove that for $\sigma \in \Stab(\A_Q)^N$.  
Let $\sigma =(Z,\sli) \in \Stab(\A_Q)^N$ and suppose $(Z,\sli) = ([\Phi]^{-1} \cdot Z,\Phi(\sli)) $. 
Since any object in $\D_Q$ is generated by finite extensions of objects given by the shift of $S_1,\dots,S_n$, 
the isomorphism 
$\Phi ( S_i) \cong S_i$ for all $i =1,\dots,n$ implies $\Phi \cong \id$. 
Assume that $S_i \in \sli(\phi_i) \subset \A_Q$, 
then also $\Phi (S_i) \in \sli(\phi_i) \subset \A_Q$ 
since $\Phi(\sli(\phi_i)) = \sli(\phi_i)$. 
At the $K$-group level, $[\Phi]^{-1} \cdot Z = Z$ implies that 
$[\Phi] = \id$, therefore we have 
$[\Phi(S_i)] = [S_i] $. Both $\Phi(S_i)$ and $S_i$ are 
objects in $\A_Q$ with the same $K$-group class $[S_i]$, and 
such an object is unique up to isomorphism in $\A_Q$. Hence $\Phi(S_i) \cong S_i$. 
  
Next we prove that the action of $\Br(\D_Q)$ is properly discontinuous. 
It is sufficient to prove that for any $\sigma \in \Stab(\A_Q)^N$, there is some open subset 
$\sigma \in U $ such that $U \cap \Phi \cdot U = \phi$ for any 
$ \id \not\cong \Phi \in \Br(\D_Q)$. 
There are two cases either $[\Phi] = \id$ or not. 
In the case $[\Phi] \neq \id$, this immediately follows from Proposition \ref{free} and the local isomorphism property of 
Theorem \ref{local_iso}. 
In the case $[\Phi] = \id$, the stability conditions $\sigma$ 
and $\Phi \cdot \sigma$ are different but
have the same central charge. 
Hence, the result follows from Lemma \ref{distance}. 
\endproof

Denote by $\underline{\pi}$ the composition of 
$\pi \colon \Stab^{\circ}(\D_Q) \to \Xreg$ and $\Xreg \to \Xreg \slash W$. 
Now we prove Theorem \ref{main}. 

{\bf Proof of Theorem \ref{main}.} The remaining part is to show  that the quotient of 
$\Stab^{\circ}(\D_Q)$ by $\Z[2] \times \Br(\D_Q)$ coincides with 
$\Xreg \slash W$. This is equivalent to that 
for any $\sigma_1,\sigma_2 \in \Stab^{\circ}(\D_Q)$, if 
$\pi(\sigma_1)=\pi(\sigma_2)$, then there are elements 
$[2n] \in \Z[2]$ and $\Phi \in \Br(\D_Q)$ with $[\Phi]=\id$ 
such that $\sigma_1 =\Phi \cdot [2n] \cdot \sigma_2$. 

By Corollary \ref{action}, we can assume that $\sigma_1 \in \Stab(\A_Q)^N$. Further, there are elements $k \in \C$ and $\Phi \in \Br(\D_Q)$ such that 
$\sigma_2^{\prime} := \Phi \cdot k \cdot \sigma_2$ lies in $\Stab(\A_Q)^N$. 
By mapping $\sigma_2^{\prime}$ onto $\Xreg$, we have 
\begin{equation*}
\pi(\sigma_2^{\prime}) = [\Phi] \cdot e^{-i \pi k} \cdot \pi(\sigma_2) 
= [\Phi] \cdot e^{-i \pi k} \cdot  \pi(\sigma_1).
\end{equation*}
Since $\Stab(\A_Q)^N$ is mapped isomorphically onto the normalized chamber  
$C^N $ and both $\pi(\sigma_1)$ and 
$\pi(\sigma_2^{\prime})$ lie in $C^N$, 
we have $[\Phi]=\id$, $k =2n \in 2\Z$ 
and $\pi(\sigma_1)=\pi(\sigma_2^{\prime})$.
\endproof

Corollary \ref{cor_aut} is proved in the completely same way as the proof of Corollary 1.5 in \cite{Br3}.


\begin{thebibliography}{CBVdB04}

\bibitem[Bri07]{Br1}
T.~Bridgeland.
\newblock Stability conditions on triangulated categories.
\newblock {\em Ann. of Math. (2)}, 166(2):317--345, 2007.

\bibitem[Bri08]{Br2}
T.~Bridgeland.
\newblock Stability conditions on {$K3$} surfaces.
\newblock {\em Duke Math. J.}, 141(2):241--291, 2008.

\bibitem[Bri09a]{Br4}
T.~Bridgeland.
\newblock Spaces of stability conditions.
\newblock In {\em Algebraic geometry---{S}eattle 2005. {P}art 1}, volume~80 of
  {\em Proc. Sympos. Pure Math.}, pages 1--21. Amer. Math. Soc., Providence,
  RI, 2009.

\bibitem[Bri09b]{Br3}
T.~Bridgeland.
\newblock Stability conditions and {K}leinian singularities.
\newblock {\em Int. Math. Res. Not. IMRN}, (21):4142--4157, 2009.

\bibitem[BS]{BrSm}
T.~Bridgeland and I.~Smith.
\newblock Quadratic differentials as stability conditions.
\newblock arXiv:1302.7030.

\bibitem[BS72]{BS}
E.~Brieskorn and K.~Saito.
\newblock Artin-{G}ruppen und {C}oxeter-{G}ruppen.
\newblock {\em Invent. Math.}, 17:245--271, 1972.

\bibitem[BT11]{BTh}
C.~Brav and H.~Thomas.
\newblock Braid groups and {K}leinian singularities.
\newblock {\em Math. Ann.}, 351(4):1005--1017, 2011.

\bibitem[CBvdB]{CBvdB}
W.~Crawley-Boevey and M.~Van~den Bergh.
\newblock Absolutely indecomposable representations and {K}ac-{M}oody {L}ie
  algebras.
\newblock {\em Invent. Math.}, 155(3):537--559, 2004.
\newblock With an appendix by Hiraku Nakajima.

\bibitem[IUU10]{IUU}
A.~Ishii, K.~Ueda, and H.~Uehara.
\newblock Stability conditions on {$A_n$}-singularities.
\newblock {\em J. Differential Geom.}, 84(1):87--126, 2010.

\bibitem[Kac78]{Kac2}
V.~G. Kac.
\newblock {I}nfinite-dimensional algebras, {D}edekind's {$\eta $}-function,
  classical {M}\"obius function and the very strange formula.
\newblock {\em Adv. in Math.}, 30(2):85--136, 1978.

\bibitem[Kac90]{Kac1}
V.~G. Kac.
\newblock {\em Infinite-dimensional {L}ie algebras}.
\newblock Cambridge University Press, Cambridge, third edition, 1990.

\bibitem[Kel08]{Kel}
B.~Keller.
\newblock Calabi-{Y}au triangulated categories.
\newblock In {\em Trends in representation theory of algebras and related
  topics}, EMS Ser. Congr. Rep., pages 467--489. Eur. Math. Soc., Z\"urich,
  2008.

\bibitem[Kin94]{King}
A.~D. King.
\newblock Moduli of representations of finite-dimensional algebras.
\newblock {\em Quart. J. Math. Oxford Ser. (2)}, 45(180):515--530, 1994.

\bibitem[KS]{KS}
M.~Kontsevich and Y.~Soibelman.
\newblock Stability structures, motivic {D}onaldson-{T}homas invariants and
  cluster transformations.
\newblock arXiv:0811.2435.

\bibitem[Lus91]{Lus}
G.~Lusztig.
\newblock Quivers, perverse sheaves, and quantized enveloping algebras.
\newblock {\em J. Amer. Math. Soc.}, 4(2):365--421, 1991.

\bibitem[Oka06]{Ok2}
S.~Okada.
\newblock On stability manifolds of {C}alabi-{Y}au surfaces.
\newblock {\em Int. Math. Res. Not.}, pages Art. ID 58743, 16 pages, 2006.

\bibitem[Par]{Pa}
L.~Paris.
\newblock $k(\pi,1)$ conjecture for artin groups.
\newblock arXiv:1211.7339.

\bibitem[ST01]{ST}
P.~Seidel and R.~P. Thomas.
\newblock Braid group actions on derived categories of coherent sheaves.
\newblock {\em Duke Math. J.}, 108(1):37--108, 2001.

\bibitem[Tho06]{T}
R.~P. Thomas.
\newblock Stability conditions and the braid group.
\newblock {\em Comm. Anal. Geom.}, 14(1):135--161, 2006.

\bibitem[vdL83]{vdL}
H.~van~der Lek.
\newblock The homotopy type of complex hyperplane complements.
\newblock {\em Ph. D. Thesis,Nijmegen}, 1983.

\end{thebibliography}
\end{document}